\documentclass[a4paper,12pt]{amsart}
\usepackage{times,a4wide,mathrsfs,amsthm,amssymb}
 \usepackage[cal=boondoxo]{mathalfa} 
\usepackage{hyperref}
\usepackage{color}

\newtheorem{thm}{Theorem}[section]

\newtheorem*{thm*}{Theorem}
\newtheorem*{corr*}{Corollary}
\newtheorem{lemma}[thm]{Lemma}

\newtheorem*{prop*}{Proposition}

\newtheorem{corr}[thm]{Corollary}
\theoremstyle{definition}

\theoremstyle{remark}

\newtheorem{rmk}[thm]{\textit{Remark}}
\renewcommand{\proof}{\noindent\textit{Proof}\/: \,\,}

\newtheorem{nonumberingt}{Acknowledgements}

\def\bZ{\mathbf{Z}}
\def\bC{\mathbf{C}}

\newcommand\LL{{\mathcal L}} 
 
\newcommand\OO{{\mathcal O}} 
\renewcommand\SS{{\mathcal S}}

\newcommand{\comp}{\raise1pt\hbox{{$\scriptscriptstyle\circ$}}}
 
\newcommand{\underint}{\int_{\raise-4pt\hbox{\hskip-8pt $-$}}}
\def\lset{\{}  
\def\rset{\}}  
\def\set#1{\lset#1\rset} 
\def\st{\mid}   
\def\sett#1#2{\lset #1 \st #2 \rset}  

\def\tr{\mathop{\rm Tr}\nolimits}
\newcommand\Tr{{}^{\mathsf{T}}\kern-0.9pt} 

  \def\mapright#1{\mathop{\vbox{\ialign{
                                ##\crcr
    ${\scriptstyle\hfil\;\;#1\;\;\hfil}$\crcr
 \noalign{\kern2pt\nointerlineskip}
    \rightarrowfill\crcr}}\;}}

\def\mapleft#1{\mathop{\vbox{\ialign{
                                ##\crcr
    ${\scriptstyle\hfil\;\;#1\;\;\hfil}$\crcr
 \noalign{\kern2pt\nointerlineskip}
    \leftarrowfill\crcr}}\;}}
\def\into{\hookrightarrow}

\def \alb{\text{Alb}}

\def\chow{\mathop{\mathsf{CH}}\nolimits}
\def\Cor{\mathop{\rm Corr}\nolimits}
\newcommand\End{\mathop{\rm End}\nolimits}
\def\half{\frac 12} 

\newcommand\im{\operatorname{Im}}
\def\jo{\iota}
\date{7 September   2017} 

\title{A motivic study of generalized Burniat surfaces}
\author[Chris Peters]
{Chris Peters}
\address{Discrete Mathematics, Technische Universiteit Eindhoven,
Postbus 513,
5600 MB Eindhoven, Netherlands.}
\email{c.a.m.peters@tue.nl}

\begin{document}
\begin{abstract}  Generalized Burniat surfaces are surfaces of general type with $p_g=q$ and Euler number $e=6$ obtained by  a variant of Inoue's 
construction method for the classical Burniat surfaces.
I prove a variant of the Bloch conjecture for these surfaces. The method applies also to the so-called Sicilian surfaces introduced in \cite{BCF}.
This implies that the Chow motives of all of these surfaces are finite--dimensional in the sense of Kimura.
\end{abstract}



\maketitle

\section{Introduction} 
Quite recently  in \cite{BCF} I. Bauer et. al. have investigated a generalization of Inoue's construction \cite{inoue} of the classical  Burniat surfaces \cite{bu}.
These surfaces are minimal, of general type and  have invariants $p_g=q=0,1,2$ or $3$.

Recall that  the Chow group $\chow_k$  is said to be  "trivial", if the natural cycle class map $\chow_k \into H_{2k}$ is injective. 
The kernel of the cycle class map $\chow_k^{\rm hom}$  can be investigated by means of the Abel--Jacobi map $\chow_k^{\rm hom} \to J^k$,
where the target is the $k$-th intermediate Jacobian. Its kernel is denoted $\chow_k^{\rm AJ}$. If this vanishes,
this has strong  consequences. For instance for surfaces this implies $p_g=0$ and the Albanese map is an isomorphism up to torsion.
The converse is Bloch's conjecture  \cite{B}.
In a follow-up study \cite{BF}  this conjecture has been verified  for the generalized Burniat surfaces, i.e  $\chow_0^{\rm hom}=0$.

These  generalized Burniat surfaces  $Y=X/G$ are all quotients of  $X$ by a freely acting abelian group $G\simeq (\bZ/2\bZ)^3$ and where
$(X,G)$ is a so called
\emph{Burniat hypersurface pair} $(X,G)$: $X$ is a  hypersurface in a product  $A$ of three elliptic curves  having at most nodes as singularities and   
$G$ is an  abelian group   acting freely on $A$ and leaving $X$ invariant.  The surface $X$ is also called \emph{Burniat hypersurface}.
These  come in $16$ families, enumerated  $\SS_1$--$\SS_{16}$.
The classical Burniat surface belongs to the $4$-parameter family $\SS_2$. Also the family $\SS_1$
is $4$-dimensional. The remaining families have only $3$ parameters coming from varying the
elliptic curve. This implies that the equation of $X$ in these cases is uniquely determined,
contrary to the first two where there is a pencil of hypersurfaces invariant under $G$. The surfaces $Y$ have at most nodal singularities.
For simplicity  I assume in this note that $Y$, and hence $X$ is smooth, which is generically the case. However, none of the arguments is influenced by
the presence of nodal singularities.

In  \cite{lnp}  it  has been remarked  that   the  main theorem of loc. cit. 
yields the Bloch conjecture for the classical Burniat surfaces.
The goal of this paper is to apply the same  methods to all  Burniat hypersurfaces.  In particular, one obtains a short proof of the Bloch conjecture
in the appropriate cases.

To state the result, let me recall that the Chow motive $h(X)$ is the pair $(X,\Delta)$ where $\Delta\subset X\times X$ is the
diagonal considered as a  (degree $0$)  self-correspondence of $X$. As a self-correspondence it is an idempotent in the ring $\Cor^0(X)$ of degree $0$
self-correspondences.

 If a finite group $G$ acts on $X$, any character $\chi$ of the group defines an idempotent
\[
\pi_\chi= \frac{1}{|G|} \sum_{g\in G} \chi(g) \Gamma_g \in \Cor^0(X),
\]
where $\Gamma_g$ is the graph of the action of $g$ on $X$. 
The pair $(X,\pi_\chi)$ is the  motive canonically associated to the character $\chi$. Note that the trivial character gives  the motive $h(X/G)$ of the
Burniat hypersurface.
The main result now reads as follows:

\begin{thm*}[=Theorem~\ref{thm:Main}] With $i:X\into A$ the inclusion, 
let $(X,G)$ be a Burniat hypersurface pair as before and let $Y=X/G$ be the corresponding generalized Burniat surface.
Consider   the one-dimensional space $H^0(\Omega^3_A)$  as a $G$-representation space, i.e. 
as  a character $\chi_A$. 
Then 
\begin{enumerate}
\item For the families $\SS_1,\SS_2$ the involution $j=\jo_1\jo_2\jo_3$ belongs to $G$ and the motives  $h(X/j)$ and   $ h(Y)$ are  finite dimensional.
 \item For all   other  families, the motive $(X,\pi_{\chi_A})$ is finite-dimensional. For the families $\SS_3$,$\SS_4$, $\SS_{11},\SS_{12},\SS_{16}$
this motive is  just $h(Y)$. 
\item The Bloch conjecture holds for the families $\SS_1$--$\SS_4$. In the remaining cases a variant of Bloch's conjecture holds,
namely\footnote{As a matter of notation, for  any $G$-module $V$ we set
$
V^\chi := \sett{ v\in V}{g (v) =\chi(g) v \text{ for all } g\in G} 
$.}   $\ker (i_* :\chow_0^{\rm AJ}(X)^{\chi_A} \to \chow_0^{\rm AJ}(A)^{\chi _A} )=0$. For the families $\SS_{10},\SS_{11}$ and $\SS_{16}$ this 
means that $\ker (i_* :\chow_0(X)  \to \chow_0(A) )=0$.\footnote{This can also be stated directly in terms of the so-called "variable motive". See
Theorem~\ref{thm:Main}.}
\end{enumerate}
\end{thm*}

As shown in \cite{BCF} the families $\SS_{11}$ and $\SS_{12}$  give two divisors in a component of the  $4$-dimensional Gieseker moduli space.

The above  method applies  also to  the surfaces in this component, the so-called \emph{Sicilian surfaces} so that 
 the result for $\SS_{11}$ and $\SS_{12}$  is valid for these as well.
See Remark~\ref{Sicilian}.

 \begin{nonumberingt} Thanks to Robert Laterveer for his  interest and  remarks. 
\end{nonumberingt}

\section{Preliminaries}

\subsection{A criterion  for finite dimensional motives}

The general situation of \cite{lnp} concerns smooth $d$-dimensional complete intersections $X$ inside a smooth projective manifold $M$
of dimension $d+r$ 
for which Lefschetz' conjecture $B(M)$ holds. This conjecture  is known to hold for projective space and for abelian varieties and so
in particular for the situation in this note.

Recall also that in this situation, with $i:X\into M$ the inclusion, the fixed and variable cohomology is defined as follows.
\[
\aligned
H^d_{\rm fix}(X) &=\im( i^* : H^{d}(M) \to H^d(X)),\\
H^d_{\rm var}(X) &= \ker(i_*:H^{d} (X)\to H^{d+2r}(M)),
\endaligned
\]
and that one has  a direct sum decomposition
\[ 
H^d(X)=
H_{\rm fix}^d(X)   \oplus H_{\rm var}^d(X),
\]
which is orthogonal with respect to the intersection product.   In \cite{MotVar} I explained that validity of $B(M)$ implies the existence
of a motive $h(X)^{\rm var}= (X,\pi^{\rm var})$ such that $\pi^{\rm var}$ induces projection onto variable cohomology.

The main input is the  special case of  \cite[Thm. 6.5 and Cor. 6.6]{lnp} for surfaces inside a threefold. It reads as follows.
\begin{thm} \label{main2Bis}   Let $M$ be a smooth projective threefold on which a finite abelian group $G$ acts.
Let $\LL$ be line bundle with $G$-action, $X\subset M$ a $G$-invariant section and $\chi$ a character of $G$.  Suppose that
\begin{enumerate}
\item  the conjecture $B(M)$ holds;
\item   the sections of $H^0(M,\LL)^G$    separates  orbits;
\item   all characters of $G$   appear in $\End(H^{2}_{\rm var}(X))$;
\item  the Chow motive $(M,\pi_\chi)$ is finite-dimensional;
\item $ 0 \not= H^2(X)_{\rm var}$ and  $H^2_{\rm var}  (X)_{\chi} \subset{ N}^1 H^2(X)$.  
\end{enumerate} 
Then    $\chow_0^{\rm var}(X)^{\chi}=0$,  and the motive  $(X,\pi_\chi)$ is finite-dimensional.   
\end{thm}

\subsection{Elliptic curves}
Let me    recall  the relevant facts  about theta functions on an elliptic curve $E$ with period lattice $\Lambda$ generated by $1$ and $\tau \in \mathfrak h$.
Points  in the elliptic curve  referred to by the standard coordinate   $z\in \bC$ and the corresponding line bundle by $\LL_z$. It is the bundle with
$H^0(E,\LL_z)= \bC \vartheta_z$, $\vartheta_z$ a  theta-function with simple  zeros in the points $z +  \Lambda$ only.  Let $t_u : z \mapsto z+u$ be a translation of $E$.
Then $ \LL_z \simeq  t_z^* \LL_0$.  If  one takes for $z$ one of the four two-torsion points  $\epsilon \in \set{0, \half,\half\tau, \half+\half\tau}$  of $E$,   the line corresponding
line bundles $\LL_\epsilon$  have  the four  classical theta functions $\vartheta_1, \vartheta_2,\vartheta_3, \vartheta_4$  respectively as sections.
See e.g. \cite[Appendix A, Table 16 ]{ency} for the definitions. 

Set
\[
M_E:= H^0(E,\LL_0^2).
\]
\begin{lemma} \label{lem:repspace}
{\rm i)} The bundle $\LL_0^2$ is a symmetric line bundle and all its sections are symmetric. \\
{\rm ii)}  The translations $t_\epsilon$ define a  faithful action action of $(\bZ/2\bZ)^2$ on $\LL_0^2$.\\
{\rm iii)} The character decomposition of $M_E$ for this action is  $(+-)\oplus  (-+)$.
\end{lemma}
\proof  i) is clear. \\
ii) Since  $\LL_\epsilon^2 \simeq  \LL_0^2$ for all two-torsion points $\epsilon$, the functions  $\vartheta_j^2$ define sections of the same bundle  
$\LL_0^2$.   The  sections  $\vartheta_j^2$, $j=1,2,3,4$ are characterized by having a double zero at  exactly one of the four  $2$-torsion points.
This shows in particular that the action of the group $\set{ t_\epsilon, \epsilon \text{ a $2$-torsion point}}$ is faithful on $M_E$.\\
iii)
  It follows that there is a  basis of   two sections of $\LL_0^2$  consisting of simultaneous eigenvectors for this action. Since the action is faithful, 
the character decomposition must be $(+-), (-+)$.
\qed\endproof


 \section{Surfaces inside abelian threefolds invariant under involutions}

\subsection{Products of three elliptic curves}
Consider the abelian threefold
\[
A:= E_1 \times E_2 \times E_3  , \quad E_\alpha= \bC /\Lambda_\alpha, \, \text{ with }   \Lambda_\alpha=\bZ\oplus \bZ\tau_\alpha,\, \alpha=1,2,3.
\]
Using for a fixed   elliptic curve   $E=\bC \bZ\oplus \tau\bZ$  the involutions
\[
\iota_E : z\mapsto -z, \quad  t_E    : z\mapsto - z+{\half} ,\quad {\tau }_E : z \mapsto  - z+ {\half\tau},
\] 
we obtain  three involutions on $A$
 \begin{equation} \label{eqn:3Invols}
 \aligned
 \iota_\alpha & =\iota_{E_\alpha}   \\
\iota_{\alpha\beta}  & = t_{E_\alpha}  t_{E_\beta},\\
\iota_{123} &=  \tau_{E_1}\tau_{E_2}\tau_{E_3}.
 \endaligned 
 \end{equation}
and we consider the group $(\bZ/2\bZ)^6$ operating on $A$ as
\[
G_0 := \langle \,   
\iota_1, \iota_2, \iota_3 ,\iota_{12} , \iota_{13}, \iota_{123}  \,  \rangle.
\]

\begin{lemma} \label{lem:ActionOnOneforms}
The action of $G_0$ on  holomorphic $1$-forms  of $A$  is given by
\begin{center}
\begin{tabular}{|c|c|c|c||c|c|c||c|c| c|}
\hline
{\rm form} &    $\iota_1$& $\iota_2$&$\iota_3$& $\iota_{12}$  &  $\iota_{13}$ &$\iota_{23}$& $\iota_{123}$ \\
\hline
$dz_1$ & $- $  & $+$    & $+$ & $-$ &  $-$& $+$   &  $- $\\
\hline
$dz_2$  & $+ $ &  $-$ & $+$ & $- $    & $+$& $ -$&   $- $ \\
\hline
$dz_3$  & $+ $ &  $+$ & $-$ & $+ $    & $-$  & $ -$&$ -$ \\
\hline
\end{tabular}
\end{center}

\end{lemma}

Consider now the symmetric line bundle   $\LL_A^2$ where
$$
\LL_A  := \OO_{E_1}( \LL_0) \boxtimes \OO_{E_2}(  \LL_0) \boxtimes  \OO_{E_3}( \LL_0), 
$$
and set 
\[
H^0(\LL_A^2) =  M_{E_1}  \boxtimes  M_{E_3} \boxtimes  M_{E_3}. 
\]
By Lemma~\ref{lem:repspace} this is a representation space for $G_0$ which admits a basis   of 
simultaneous eigenvectors. If \set{$\theta^1_{E_j}, \theta^2_{E_j}}$,  denotes the basis  of Lemma~\ref{lem:repspace}, for $M_{E_j}$, $j=1,2,3$, 
their  products give  $8$ basis vectors as follows.
\[
\theta_{j_1j_2j_3}= \theta_{E_1}^{j_1} \cdot  \theta_{E_2}^{j_2} \cdot \theta_{E_3}^{j_3},\quad j_k\in{1,2}.
\]
The next result is a consequence.:
\begin{lemma} \label{lem:GactionOnTheta}
The space $H^0(\LL^2_A)$ is the $G_0$-representation space which on the basis $\{ \theta_{j_1j_2j_3}\}$,  $j_1,j_2,j_3\in \{1,2\}$, 
is given as follows
\begin{center}
\begin{tabular}{|c|c|c|c||c|c|c||c|}
\hline
{\rm element} &   $\iota_1$& $\iota_2$&$\iota_3$& $\iota_{12}$  &  $\iota_{13}$ &$\iota_{23}$& $\iota_{123}$ 
\\
\hline
$\theta_{111}$ & $+ $ &  $+$   & $+$ & $+ $ & $+$ &    $+$ &   $-$    \\
\hline
\hline
$\theta_{211}$  & $+ $ &  $+$   & $+$       & $- $ & $-$   &  $+$  & $ +$     \\
\hline
$\theta_{121}$   & $+ $ &  $+$ &   $+$ &    $-$  & $+$  &  $-$  & $ +$    \\
\hline
$\theta_{112}$ & $+ $ &  $+$  &    $+$ &    $+ $ & $-$ &   $-$   &$ +$    \\
\hline
\hline
$\theta_{1\, 2\, 2}$  & $+ $ &  $+$ & $+$ & $- $ & $-$  & $+$   & $-$   \\
\hline
$\theta_{2\, 1\, 2 }$   & $+ $ &  $+$ & $+$ & $- $  & $+$ & $-$   & $-$   \\
\hline
$\theta_{221}$ & $+ $ &  $ +$ & $+$ & $+$ & $-$&   $-$  & $-$ \\
\hline
\hline
$\theta_{222}$  & $+ $ &  $+$  & $+$ & $+$ & $+$&  $+$   & $ +$   \\
\hline
\end{tabular}
\end{center}
\end{lemma}

\subsection{Hypersurfaces of abelian threefolds and involutions}

Let $A$ be an abelian variety of dimension three and $\LL$ a principal polarization so that $\LL^3=3! =6 $ and
let $i: X\into  A$ be a smooth surface given by a section of of $\LL^{\otimes 2}$.
The Lefschetz's hyperplane theorem  gives:
\begin{eqnarray}
i^* : H^1(A)& \mapright{\simeq} &H^1(X  \label{eqn:Lef1})\\
i^*:H^2(A)& \mapright{\simeq} &H^2_{\rm fix} (X)\subset H^2(X) \label{eqn:Lef2}.
\end{eqnarray}

\begin{lemma} \label{lem:TrOnHVar}
Suppose that $\iota :A\to A$ is an involution which acts on $H^0(\Omega^1_A)$ with $p$ eigenvalues $1$ and $n=3-p$ eigenvalues $-1$. Suppose also 
 that $\iota$ preserves $X$ and acts without fixed points on $X$. Then we have 
 \[
\tr(\iota) |H^2_{\rm var}(X) = -29 +8p(4-p) = \begin{cases} -29 & \text{for }  p=0\\
 										 -5 & \text{for }  p=1 \\
										  \,\,3  & \text{for }  p=2  \\
										  \, -5 &  \text{for }  p=3.  \end{cases}
\]
 \end{lemma}
\proof The assumption implies that
\[
\tr(\iota)|H^1(A) = 4p - 6= \begin{cases} -6  & \text{for }  p=0\\ 
							     -2 & \text{for }  p=1 \\
							      \, 2 &\text{for }  p=2\\ 
							     \, 6& \text{for }  p=3, \end{cases}
							    \]
							    and
							    \[
 \tr(\iota)| H^2(A)= 8p(p-3)+15=  \begin{cases} 15 & \text{for }  p=0\\
  -1 & \text{for }  p=1 \\ 
  \, 1& \text{for }  p=2 \\
  \, 15& \text{for }  p=3. \end{cases}
\]
If $\iota$ preserves $X$ and acts without fixed points on $X$, Lefschetz' fixed point theorem gives
\[
0= 2 -  2 \tr(\iota) | H^1(X) + \tr(\iota) | H^2(X)=  2 -  2\tr(\iota) | H^1(A) + \tr(\iota) | H^2(A) + \tr(\iota) |H^2_{\rm var}(X) ,
\]
and so the above calculation immediately gives the desired result. \qed\endproof

In order to calculate the invariants on $X$, let me first  consider the holomorphic two-forms in detail.
\begin{lemma}  \label{lem:ActionOnHolForms} {\rm 1.} One has
\[
h^{0,2}_{\rm var} (X) = 7,\quad h^{0,2}_{\rm fix} (X)= 3.
\]
{\rm 2.}  If $X=\{\theta_0=0\}$,  the variable holomorphic $2$-forms are the Poincar\'e-residues along $X$ of 
the meromorphic $3$-forms  on $A$  given by expressions of the form
\[
\frac{\theta}{\theta_0} dz_1\wedge dz_2\wedge dz_3
\]
with $\theta$  a  theta-function  on $A$ corresponding to a section  of  $\LL^{\otimes 2}$, and where  $z_1,z_2,z_3$  are holomorphic coordinates on $\bC^3$.\\
{\rm 3.} Suppose $\jo$ acts with the  character $\epsilon\in \set{\pm 1}$ on holomorphic three forms.  Let $(p,n)$ be the dimensions of the invariant, resp. anti-invariant
sections of $\LL^{\otimes 2}$. Then  $\dim H^{2,0}_{\rm var,+}(X)= p-1$ if $  \epsilon=1 $ and $=p$ otherwise.
\end{lemma}
\proof
Consider  the Poincar\'e residue sequence
\[
0 \to \Omega^3_A  \to \Omega^3_A(X)  \mapright{\rm res}  \Omega^2_X \to 0.
\]
In cohomology this gives
\begin{equation}\label{eqn:PR}
0 \to H^0(\Omega^3_A) \to H^0(\Omega^3_A(X)) \mapright{\rm res}  H^0(\Omega^2_X) \to H^1(\Omega^3_A) \to 0.
\end{equation}
Since $H^0(\Omega^3_A(X))= H^0(\LL^{\otimes 2})$ the  assertion 1. follows.\\
2. This is clear. \\
3.  This follows directly from \eqref{eqn:PR}. \qed \endproof

\begin{corr}
The invariants of $X$ are as follows.
\begin{center}
\begin{tabular}{|c|c|c|}
\hline
 $b_1$ & $b_2^{\rm var}=(h^{2,0}_{\rm var},h^{1,1} _{\rm var},h^{0,2} _{\rm var} )$  & $ b_2^{\rm fix}=(h^{2,0}_{\rm fix},h^{1,1} _{\rm fix},h^{0,2} _{\rm fix} ) $ \\
\hline
   $6$ & $43=(7,29,7)$ & $15=(3,9,3)$ \\
\hline
\end{tabular}
\end{center}
\end{corr}
\proof  Equation~\eqref{eqn:Lef1}  gives $b_1(X)=b_1(A)=6$. 
To calculate $b_2(X)$ we observe that $c_1(X)= -2 \LL|_X$ and $c_2(X)=  4 \LL^2|_X$ so that 
\[
c_1^2(X)=c_2(X)=4 \LL^2|_X=8 \LL^3 =48.
\]
Since    $c_2(X)= e(X)=2- 2b_1(X)+ b_2(X)=48$, it follows that $b_2(X)= 58$.  By \eqref{eqn:Lef2} one has   $b_2^{\rm fix,+}(X)= b_2(A)=15$ and so $b_2^{\rm var}(X)= 43$.
Since $h^{2,0}_{\rm var}=7$, the invariants for $X$ follow.
\qed\endproof

\subsection{Burniat hypersurfaces}

A   \emph{Burniat hypersurface}  of $A$ is a surface which is invariant under a subgroup $G\subset G_0$ generated by $3$ commuting involutions and which acts freely
on $X$. Each of the involutions is a product of  the involutions  \eqref{eqn:3Invols}.
The quotient $Y=X/G$ is called a generalized Burniat surface.
In \cite{BCF} one finds a list of $16$ types of such surfaces, denoted $\SS_1,\dots,\SS_{16}$.  All of the surfaces are of general type with $c_1^2=6$, $p_g=q$ and $q=0,1,2,3$
and hence $e=6= 2 -4q +b_2$ so that $b_2=(p_g,h^{11},p_g)= (q, 4+2q, q)$.
There are $4$ families with $q=0$ and one of these, $\SS_2$   gives  the classical Burniat surfaces from \cite{bu}. See  Table~\ref{table:BurniatHyps}. 
\begin{table}[htp]
\caption{Burniat hypersurfaces} \label{table:BurniatHyps} 
\begin{center}
\begin{tabular}{|c|c||c||c||c||c|}
\hline
type &   involution $1$  & involution $2$  & involution $3$ &$G$-invariant $1$-forms& $\chi_A$ \\
\hline
$\SS_1$ & $\iota_1\iota_2\iota_3$  & $\iota_2\iota_3\iota_{123}$ & $\iota_3\iota_{23}$   &none& $---$\\ 
\hline
$\SS_2$ & $\iota_1  \iota_3   \iota_{23}$ &  $\iota_3\iota_{13}$  & $\iota_2\iota_{23}$    & none& $+--$  \\
\hline
$\SS_3$  & $\iota_1\iota_3\jo_{23}$& $\iota_3\iota_{123}$ & $\iota_2\iota_3\iota_{12}$ & none & $+++$ \\
\hline
$\SS_4$ & $\iota_1\iota_3\iota_{12}$ & $\iota_2\iota_{123}$ & $\iota_2\iota_3\iota_{23}$ & none &   $+++$ \\
\hline
\hline
$\SS_5$ & $\jo_1\jo_3\jo_{13}$ & $\jo_3\jo_{123}$ & $\jo_3\jo_{23}$ & $dz_3$  & $++-$  \\
\hline
$\SS_6$  & $\jo_2\jo_3\jo_{123}$ & $\jo_2\jo_3\jo_{13}$ & $\jo_3\jo_{23}$ & $dz_3$ &   $-+-$\\
\hline  
$\SS_7$  & $\jo_1\jo_3\jo_{23}   $ &  $\jo_3\jo_{123}$ & $\jo_2\jo_{12}  $ & $dz_3$ &$++- $  \\
\hline
$\SS_8$ & $\jo_1\jo_3\jo_{23}$ & $\jo_2\jo_3\jo_{123}$ & $\jo_2\jo_3\jo_{13}$ &$dz_3$ &  $+-+$  \\
\hline  
$\SS_9$  & $\jo_1\jo_2\jo_3\jo_{13} $ & $\jo_3\jo_{123}$ & $\jo_2\jo_{12}$ &  $dz_3$ &  $-+-$ \\
\hline 
$\SS_{10}$  & $\jo_1\jo_2\jo_3\jo_{13}$ & $\jo_2\jo_3\jo_{123}$ & $\jo_3\jo_{23}$ &$dz_3$ & $---$ \\
\hline   
$\SS_{11}$ & $\jo_1\jo_2\jo_{23} $ & $\jo_2\jo_{123}$ & $\jo_2\jo_3\jo_{12}$ & $dz_2$&  $+++$ \\
\hline  
$\SS_{12}$ & $\jo_1\jo_3\jo_{13} $ & $\jo_3\jo_{123}$ & $\jo_2\jo_3\jo_{23}$ &$dz_3$ &  $+++$ \\
\hline 
\hline 
$\SS_{13}$  &  $\jo_1\jo_2\jo_3\jo_{23}$ & $\jo_2\jo_3\jo_{123} $ &  $\jo_2\jo_{12}$ &  $dz_2,dz_3$ & $---$ \\
\hline
$\SS_{14}$ & $\jo_1\jo_{13} $ & $\jo_{12}\jo_{123}$ & $\jo_2\jo_{23}$ & $dz_1,dz_2$ & $---$ \\
\hline  
$\SS_{15}$ & $ \jo_1\jo_3\jo_{13} $ & $\jo_{12}\jo_{123}$ & $\jo_2\jo_3\jo_{23}$ & $dz_1, dz_2$ & $+-+$ \\
\hline
\hline  
$\SS_{16}$  & $\jo_1\jo_3\jo_{13}$   &  $\jo_3\jo_{12} \jo_{123}$ & $\jo_2\jo_3\jo_{23}$& all & $+++$ \\
\hline
\end{tabular}
\end{center}
\end{table}
The last column of this table gives the action of the three generators $(g_1,g_2,g_3)$ on $H^0(\Omega^3_A)$.
It is calculated using Lemma~\ref{lem:ActionOnOneforms}.
\begin{table}[htp]
\caption{Action on forms  and invariants of the generalized Burniat surfaces} \label{table:ActOnForms}
 \begin{center}
\begin{tabular}{|c|c||c||c|c|c|}
\hline
type &  $U=H^0(A,\Omega^1_A)$&  $H^0(A,\Omega^2_A)=\wedge^2 U=W$   &   $H^1(A,\Omega^1_A)$  & $b_2^{\rm fix}(Y)$ &  $b_2^{\rm var}(Y)$ \\
\hline
$\SS_1$ & $(-++)(--+)(-+-)   $   &$(+-+)(++-)(+-- )$ &  $ 3\cdot \mathbf 1 + 2 W$  &  $(0,3,0) $ &$(0,4,0)$ \\
\hline
$\SS_2$ & $(-++)(+-+)(++-) $   &$(+--)(-+-) (--+)$ &  $ 3\cdot \mathbf 1 + 2 W$  &  $(0,3,0) $ &$(0,4,0)$ \\
\hline
$\SS_3$ & $(---)(--+)(++-) $   &$(++-)(--+) (---)$ &  $ 3\cdot \mathbf 1 + 2 W$  &  $(0,3,0) $&$(0,4,0)$  \\
\hline
$\SS_4$ & $(+-+)(-++)(--+) $   &$(--+)(-++) (+-+)$ &  $ 3\cdot \mathbf 1 + 2 W$  &  $(0,3,0) $&$(0,4,0)$  \\
\hline
\hline
$\SS_5$ & $(+-+)(+--) +\mathbf 1 $   &$(+-+)(+--) (++-)$ &  $ 3\cdot \mathbf 1 + 2 W$  &   $(0,3,0)$&$(1,3,1)$  \\
\hline
$\SS_6$ & $(--+)(+--) +\mathbf 1 $   &$(--+)(+--)  (-+-)$ &  $ 3\cdot \mathbf 1 + 2 W$  &   $(0,3,0)$ &$(1,3,1)$ \\
\hline
$\SS_7$  & $(---) (-++) +\mathbf 1$  &  $ (---)(-++)   (+--)$ &   $ 3\cdot \mathbf 1 + 2 W$  & $(0,3,0)$&$(1,3,1)$  \\
\hline
$\SS_8$ & $(---)(-++)  +\mathbf 1 $   &$(---)(-++)  (+--)$ &  $ 3\cdot \mathbf 1 + 2 W$  &   $(0,3,0)$&$(1,3,1)$  \\
\hline
$\SS_9$ & $(+--)(--+) +\mathbf 1 $   &$(+--)(--+)(-++)$ &  $ 3\cdot \mathbf 1 + 2 W$  &   $(0,3,0)$&$(1,3,1)$  \\
\hline
$\SS_{10}$ & $(+-+)(-+-)+\mathbf 1  $   &$(+-+)(-+-) (---)$ &  $ 3\cdot \mathbf 1 + 2 W$  &  $(0,3,0)$ &$(1,3,1)$ \\
\hline
$\SS_{11}$  & $2(---)   +\mathbf 1$  &  $ 2(---)+\mathbf 1 $ &   $ 3\cdot \mathbf 1 + 2 W$  &  $(1,5,1)$ &$(0,1,0)$  \\
\hline
$\SS_{12}$  & $2(+-+) +\mathbf 1$  &  $2(+-+) +\mathbf 1$ &   $ 3\cdot \mathbf 1 + 2 W$  &  $(1,5,1)$ &$(0,1,0)$  \\
\hline
\hline
$\SS_{13}$  & $(+--)+ 2 \cdot \mathbf 1$     &$2 \cdot (+--) +\mathbf 1   $ & $   3\cdot \mathbf 1 + 2 W $&  $(1,5,1)$ &$(1,3,1)$    \\
\hline
$\SS_{14}$  & $(---)+ 2 \cdot \mathbf 1$     &$2 \cdot (---) +\mathbf 1  $&   $ 3\cdot \mathbf 1 + 2 W$  &  $(1,5,1)$ &$(1,3,1)$ \\
\hline
$\SS_{15}$  & $(+-+)+ 2 \cdot \mathbf 1$     &$2 \cdot +-+) +\mathbf 1  $ &   $ 3\cdot \mathbf 1 + 2 W$  &  $(1,5,1)$ &$(1,3,1)$ \\
\hline
\hline
$\SS_{16}$  &$  3\cdot \mathbf  1$  &   $ 3\cdot \mathbf  1  $   &   $ 3\cdot \mathbf 1 + 2 W$  &   $(3,9,3$&$(0,1,0)$  \\
\hline
\end{tabular}
\end{center}
\end{table}

In Table~\ref{table:ActOnForms}  the character spaces for the action on the forms coming from $A$ is given. It is calculated from
the description of the generating involutions as given in Table~\ref{table:BurniatHyps} and  the known action of $1$-forms as given in Lemma~\ref{lem:ActionOnOneforms}.
From the first column of Table~\ref{table:ActOnForms} one finds the trace of the action of these generators on $H^0(\Omega^1_A)$, or, alternatively,
the dimensions of the eigenspaces for the eigenvalues $+1$ and $-1$.
 Writing  for example the dimensions of the $(+)$-eigenspaces as a vector according to 
  the   group elements  written in the order    $(1, g_1,g_2,g_3,g_1g_2,g_1g_3,g_2g_3,g_1g_2g_3)$ yields 
 the \emph{type} $(3,t_1,t_2,t_3,\dots)\in \bZ^8$ of the group action. 
 This   gives the  first row in Table \ref{table:TraceOnForms} below.
\begin{table}[htp]
\caption{Trace vectors} \label{table:TraceOnForms}
 \begin{center}
\begin{tabular}{|c|c||c||c|c|}
\hline
type &  type $H^0(A,\Omega^1_A)$ &   trace vect. $H^2_{\rm var}(X)$   &   trace  vect. $H^0(A,\Omega^3_A)$ &mult. $\chi_A$  \\
\hline
$\SS_5$ & $(3| 3\, 1\,1| 1\,2\,2| 2)$   &$(43|{-5\,} \,{-5\,}\,{-5\,} | {-5\,} 3\, 3| 3)$ &  $ (1|1\, 1{-1\,}|1\,{-1\,} {-1\,}|{-1)}$ & $3$   \\
\hline
$\SS_6$ & $(3| 2 \,1 \, 2 |  2  \,1 2| 2)$   &$(43  |3\,{-5\,}3 |  \,{-5\,} 3\,{-5\,} |3)$ &  $ (1|{-1\,}1\,{-1\,}  |  {-1\,}1\, {-1}|1)$  & $6$   \\
\hline
$\SS_7$  & $(3| 1\, 2\, 2|2\, 2\, 3 |1)$  &  $ ( 43 | \,{-5\,}33| 33\,{-5\,}  |  \,{-5\,}  ) )$ &   $(1|1\, 1{-1\, }|1\, {-1}{-1}|{-1}) $  & $6$  \\
\hline
$\SS_8$ & $(3| 1\, 2\, 2|2\, 2\, 3 |1)$    &  $( 43 | \,{-5\,}33| 33\,{-5\,}  |  \,{-5\,}  ) )$  &   $ (1|1\, {-1}\, 1|{-1}\, 1\, {-1}|{-1}) $   & $6$   \\
\hline
$\SS_9$ & $(3|2\, 1\, 2 | 2\, 1\, 2 |3)  $   &$(43  |3\,{-5\,}3 |  \,{-5\,} 3\,{-5\,} |\,{-5\,})$  &  $ (1|{-1\,}1\,{-1}  |  {-1}\, 1\, {-1}|1)$   &  $5$ \\
\hline
$\SS_{10}$ & $(3|2\, 2\, 2|1\, 3\, 1|2)$   &$(43| \,{-5\,}   \,{-5\,}  \,{-5\,} |  \,{-5\,}  \,{-5\,}  \,{-5\,}|3))$ &  $ (1|{-1\, }{-1\, }{-1}|   1 \,1 \, 1|{-1})$ &   $5$ \\
\hline
\hline
$\SS_{13}$  &  $(3| 3\,2\, 2 | 2\, 2\, 3  | 3) $    &$(43| \,{-5\,} 3\,3 | 3\,3\,{-5\,}| \,{-5\,})  $ &  $ (1|{-1}\,{-1}\,{-1}|   1\, 1 \, 1|{-1})$  &   $6 $     \\
\hline
$\SS_{14}$  & $(3|2\, 2\, 2|3\, 3\, 3|2)$     &$(43|3\, 3\, 3|\,{-5\,}\,{-5\,}\,{-5\,}|3)$&   $ (1|{-1}\, {-1}\, {-1}|   1 \,1 \,1|{-1})$   &   $ 2$    \\
\hline
$\SS_{15}$  & $(3|3\, 2\, 3| 2\, 3\, 2|2)$     &$(43|\,{-5\,}3\,{-5\,}|3\,{-5\,}3|3)$ &   $ (1|1\,{-1\,}1|{-1}1\,{-1}|{-1}) $  &   $2$    \\
\hline
\end{tabular}
\end{center}
\end{table}
Using Lemma~\ref{lem:TrOnHVar}, this table enables  to find the multiplicity of $\chi_A$ in  $H^2_{\rm var}(X)$. 
\begin{lemma}  For each of the families $\SS_3$--$\SS_{16}$ the multiplicity of  the character  $\chi_A$  inside  $H^2_{\rm var}(X)$   is
given in the last column of Table~\ref{table:TraceOnForms} and in particular, is  non-zero.\footnote{This is also true for the two remaining families, but this will not be used.}
\end{lemma}
\proof For  each of the   families $\SS_3,\SS_4, \SS_{11}, \SS_{12}$ and $\SS_{16}$ one has $H^0(A,\Omega^3_A)=(+++)$ and $H^{1,1}_{\rm var}(Y)= H^{1,1}_{\rm var, +++}(X) =3$ or $=1$,
as one sees from Table~\ref{table:ActOnForms}. 

For the other families we argue as follows. In each case, $g\in G$, $g\not=1 $   act freely on $X$ and so one can apply Lemma~\ref{lem:TrOnHVar}   to find $\tr g|H^2_{\rm var}(X)$, given the dimension $p(g)$ of the $(+1)$-eigenspace of $H^0(\Omega^1_A)$.  This type is given in Table~\ref{table:TraceOnForms}.
The next column gives the corresponding trace vector. Then follows the  trace vector of $\chi_A$.
Now apply the  trace formula for the multiplicity of an irreducible representation inside a given
representation (see e.g. \cite[\S 2.3]{lpfg}). Let me do this  explicitly for  the family $\SS_5$.
The trace vector for $H^2_{\rm var}(X)$ is  $(43,-5,-5,-5,-5,3,3,3)$, the first number being $\dim H^2_{\rm var}(X)$. 
The representation $\chi_A==(++-)$
has  trace vector $(1,1,1,-1,1,-1,-1,-1)$ and the trace formula  reads 
\[
\frac 18 (43-5-5+5-5-3-3-3)= 3. \hfill \qed 
\]
\endproof

 \section{The main result}
 In this section I shall show that the main theorem~\ref{thm:Main} below   follows upon application of Theorem~\ref{main2Bis}. 
 First an auxiliary  result.
 \begin{lemma} \label{lem:AllCharsAppear}
 Consider for each of the families $\SS_1$--$\SS_{16}$ the space of theta-functions $H^0(\LL^2_A) $ as $G$-representation space.
This $8$-dimensional space is the direct sum for all $8$ characters of $G$ except for the families $\SS_1$ and $\SS_2$.
For these families we have
\begin{itemize}
\item for $\SS_1$ we have $H^0(L^2_A) = 2\left(  (+++)+ (++-)+ (+-+)+(+--) \right)$,
\item  for $\SS_2$ we have $H^0(L^2_A) = 2\left(  (+++)+ (+-+)+ (-+-) + (---) \right)$.
\end{itemize}
\end{lemma}
\proof This follows  from the $G$-action on  the basis $\theta_{j_1,j_2j_3}$ for  $H^0(\LL^2_A) $ which can be deduced from Lemma~\ref{lem:GactionOnTheta}.
I shall work this out for two cases: the family $\SS_2$, and for the family $\SS_6$. 
For $\SS_2$ we have  $g_1=\iota_1  \iota_3   \iota_{23}$, $g_2=\iota_3\iota_{13}$ and  $g_3=\iota_2\iota_{23}$
and for $\SS_6$ we have  $g_1=\jo_2\jo_3\jo_{123}$, $g_2=\jo_2\jo_3\jo_{13}$  and  $g_3=\jo_3\jo_{23}$,
and the action of these involutions is given  in the following table.
\begin{center}
\begin{tabular}{|c|c|c|c||c|c|c|}
\hline
{\rm element} &   $g_1=\jo_1  \jo_3   \iota_{23}$ & $g_2=\jo_3\jo_{13}$ &$g_3=\iota_2\iota_{23}$ & $g_1=\jo_2\jo_3\jo_{123}$  &  $g_2=\jo_2\jo_3\jo_{13}$ &  $g_3=\jo_3\jo_{23}$ 
\\
\hline
$\theta_{111}$ & $+ $ &  $+$   & $+$ & $- $ & $+$ &    $+$    \\
\hline
\hline
$\theta_{211}$  & $+ $ &  $-$   & $+$       & $+ $ & $-$   &  $+$        \\
\hline
$\theta_{121}$   & $- $ &  $+$ &   $-$ &    $+$  & $+$  &  $-$       \\
\hline
$\theta_{112}$ & $- $ &  $-$  &    $-$ &    $- $ & $-$ &   $+$       \\
\hline
\hline
$\theta_{1\, 2\, 2}$  & $+ $ &  $-$   & $+$   & $- $ & $+$  & $-$       \\
\hline
$\theta_{2\, 1\, 2 }$   & $- $ &  $+$ & $-$ & $- $  & $-$ & $+$       \\
\hline
$\theta_{221}$ & $- $ &  $ -$ & $-$ & $-$ & $-$&   $-$   \\
\hline
\hline
$\theta_{222}$  & $+ $ &  $+$  & $+$ & $+$ & $+$&  $+$      \\
\hline
\end{tabular}
\end{center}
\qed\endproof
 
 \begin{thm}  \label{thm:Main}
Let $(X,G)$ be a Burniat hypersurface pair as before and let $Y=X/G$ be the corresponding generalized Burniat surface.
Consider   the one-dimensional space $H^0(\Omega^3_A)$  as a $G$-representation space, i.e. 
as  a character $\chi_A$. 
Then 
\begin{enumerate}
\item For the families $\SS_1,\SS_2$ the involution $j=\jo_1\jo_2\jo_3$ belongs to $G$ and the motives  $h(X/j)$ and   $ h(Y)$ are  finite dimensional.
 \item For all   other  families, the motive $(X,\pi_{\chi_A})$ is finite-dimensional. For the families $\SS_3$,$\SS_4$, $\SS_{11},\SS_{12},\SS_{16}$
this motive is  just $h(Y)$. 
\item The Bloch conjecture holds for the families $\SS_1$--$\SS_4$. In the remaining cases a variant of Bloch's conjecture holds,
namely  $\chow_0^{\rm var}(X)^{\chi_A}=0$. For the families $\SS_{10},\SS_{11}$ and $\SS_{16}$ this means that $\chow_0^{\rm var}(X)=0$.
\end{enumerate}
\end{thm}
\proof
(1) For the family $\SS_2$ this is \cite[Example 7.3]{lnp}. The same proof goes through for the family $\SS_1$.
\\
(2) The conditions (1), (2) and (4) of Theorem~\ref{main2Bis} are verified.
Condition (5)
is a consequence of Lemma~\ref{lem:ActionOnHolForms}.3. Indeed, if all characters in  the $8$-dimensional space $H^0(A,\LL^2)$ appear once, this
result  implies that there is one character
missing in $H^{0,2}_{\rm var}(X)$,  namely the character $\chi_A$ for the holomorphic three-forms on $A$. So for this character space 
condition (5) holds.
As to (3), Lemma~\ref{lem:AllCharsAppear} states that in this case all characters appear in $H^2_{\rm var}(X)$ except maybe this missing character
$\chi_A$. But its multiplicity has been calculated in Table~\ref{table:TraceOnForms}. It is non-zero 
  and  so  condition (3) holds as well.\\
(3) This is one of the assertions of  Theorem~\ref{main2Bis}.
\qed\endproof

\begin{rmk}\label{Sicilian}
Recall the  following definition from  \cite{BCF}: a {\em Sicilian surface\/} is a minimal surface $S$ of general type with numerical invariants
 $p_g(S)=q(S)=1, c_1^2(S)=6$ for which, in addition, there exists an unramified double cover $\hat{S}\to S$ with $q(\hat{S})=3$, and such that the Albanese morphism $\hat{\alpha}\colon \hat{S}\to  \alb(\hat{S})$ is birational to its image $Z$, a divisor in its Albanese variety  with $Z^3=12$.
In loc. cit. one finds the following explicit construction.
Let $T=\bC^2/\Lambda_2$, $\Lambda_2=\bZ^2\oplus \tau_1\bZ\oplus \bZ\tau_2$
 be an Abelian surface with a $(1,2)$-polarization $\LL_2$ and let $E=\bC/\Lambda$, $\Lambda= \bZ\oplus \tau \bZ$  be an elliptic curve. Consider 
the  sections of the line bundle $\LL=\LL_0\boxtimes \LL_2$ on $A:= E\times T$  that are invariant under the action of the bi-cyclic group $K$
generated by $(e, a) \mapsto (e+ \half\tau, -a+\half\tau_1)$ and $(e, a) \mapsto (e+ \half , a+\half\tau_2)$.
These sections define hypersurfaces   $X\subset A$ and the
quotient $Y=X/K$ is a Sicilian surface and all such surfaces are obtained in this way.

Let me consider the invariants. Note that    $p_g(Y)=q(Y)=1, c_1^2(Y)=6$  implies that $h^{1,1}(Y)=h^{1,1}(X)_{++}=6$. 
 In the same manner as   for the families $\SS_{11}$ and $\SS_{12}$ one shows that
$h^{2,0}_{\rm var, ++}=0$ so that $N^1H^2(X)_{\rm var}=H^2(X)_{\rm var}$. Moreover,  likewise $h^{1,1}(X)_{\rm var,++}=1$.
In the course of the proof of  \cite[Theorem 6.1]{BCF} it is remarked that
$H^0(A, \LL)= (++)\oplus (+-)\oplus(-+)\oplus (--)$.
Clearly, $H^0(\Omega^3_A)$ is invariant under  $K$ and the residue calculus (cf. Lemma~\ref{lem:ActionOnHolForms}) shows that $H^{2,0}_{\rm var}=(+-)\oplus(-+)\oplus (--) $ and so the "missing character" $\chi_A$ is the trivial character.  Since $h^{1,1}(X)_{\rm var, ++}=1$ this missing character is present in $H^2(X)_{\rm var}$ and by Theorem~\ref{main2Bis} it follows that
for Sicilian surfaces $Y$, one has $\chow_0^{\rm var}(Y)=0$ and the motive $h(Y)$ is finite-dimensional.

\end{rmk}

\end{document}